\documentclass{amsart}
\usepackage{amssymb,euscript,amsmath, mathrsfs, amscd, mathabx}
\usepackage[dvips]{graphicx}
\usepackage[dvips]{color}
\usepackage{bm}

\newcounter{ENUM}
\newcommand{\itm}{\item}
\newenvironment{ilist}[1][0]{\renewcommand{\theENUM}{\roman{ENUM}}\renewcommand{\itm}{\addtocounter{ENUM}{1}\item[(\theENUM)]}\begin{itemize}\setcounter{ENUM}{#1}}{\end{itemize}}
\newenvironment{Ilist}{\renewcommand{\theENUM}{\Roman{ENUM}}\renewcommand{\itm}{\addtocounter{ENUM}{1}\item[(\theENUM)]}\begin{itemize}\setcounter{ENUM}{0}}{\end{itemize}}

\newcommand{\margh}[1]{}

\input xy
\xyoption{all}
\CompileMatrices

\def\ZZ{{\mathbb Z}}

\def\bn{{\bm{n}}}

\def\sE{{\mathscr E}}
\def\sF{{\mathscr F}}

\def\sL{{\mathscr L}}

\def\sO{{\mathscr O}}

\def\sW{{\mathscr W}}

\def\fg{{\mathfrak g}}

\def\Sym{\operatorname{Sym}}

\def\Spec{\operatorname{Spec}}

\def\id{\operatorname{id}}

\def\rk{\operatorname{rk}}
\def\im{\operatorname{im}}

\def\spn{\operatorname{span}}

\def\init{\operatorname{in}}

\newtheorem{thm}{Theorem}[section]

\newtheorem{lem}[thm]{Lemma}
\newtheorem{cor}[thm]{Corollary}

\theoremstyle{definition}
\newtheorem{defn}[thm]{Definition}

\newtheorem{ex}[thm]{Example}
\newtheorem{sit}[thm]{Situation}
\theoremstyle{remark}

\newtheorem{rem}[thm]{Remark}

\numberwithin{equation}{section}
\numberwithin{figure}{section}

\begin{document}
\title{Linked determinantal loci and limit linear series}
\author{John Murray and Brian Osserman}
\begin{abstract} 
We study (a generalization of) the notion of linked determinantal
loci recently introduced by the second author, showing that as with
classical determinantal loci, they are Cohen-Macaulay whenever they have
the expected codimension. We apply this to prove Cohen-Macaulayness and
flatness for moduli spaces of limit linear series, and to prove a
comparison result between the scheme structures of Eisenbud-Harris limit 
linear series and the spaces of limit linear series recently constructed
by the second author. This comparison result is crucial in order to 
study the geometry of Brill-Noether loci via degenerations.
\end{abstract}

\thanks{The second author was partially supported by 
Simons Foundation grant \#279151 during the preparation of this work.}
\maketitle

\section{Introduction}

The theory of limit linear series for curves of compact type was developed
by Eisenbud and Harris in a series of papers in the 1980's, with the
foundational definitions and results appearing in \cite{e-h1}. They were 
able to give spectacular applications (see for instance \cite{e-h5} and 
\cite{e-h6}), despite the fact that the moduli space of limit linear series
they constructed for families of curves was not proper. However, for finer
analyses, it becomes important to have a proper moduli space. This arises
for instance when one wants to carry out intersection theory calculations
on moduli spaces of linear series, as in Khosla \cite{kh1}, or when one
wants to study the geometry of moduli spaces of linear series, as in
the current work of Castorena-Lopez-Teixidor \cite{c-l-t1} and 
Chan-Lopez-Pflueger-Teixidor \cite{c-l-p-t1}. A major step in this
direction was accomplished in \cite{os20}, when an equivalent definition
of limit linear series was introduced, leading to the first proper moduli
spaces in families. However, while the new definition was shown to agree
with the Eisenbud-Harris definition on a set-theoretic level, the scheme
structures are difficult to compare directly. The Eisenbud-Harris scheme
structure is more amenable to explicit calculation, but without knowing
that the two scheme structures agree, one cannot carry through the
arguments of \cite{c-l-t1} and \cite{c-l-p-t1}, because different 
non-reduced structures in the special fiber will typically affect
more elementary aspects of the generic fiber, such as connectedness or
the genus.

In the present paper, we address this issue by showing that under the
typical circumstances considered in limit linear series arguments, the
two scheme structures do in fact agree. We further show that when they
have the expected dimension, limit linear series spaces are Cohen-Macaulay
and flat. Our arguments center around an analysis of the 
``linked determinantal loci'' introduced in Appendix A of \cite{os25} in 
order to prove smoothing theorems for limit linear series. In fact, we
study a more general definition than the one considered in \cite{os25},
which is also a generalization of classical determinantal loci. A 
preliminary definition is the following.

\begin{defn}\label{defn:s-linked} Let $S$ be a scheme, and $d,n$ be positive
integers. Suppose that $\sE_1,\dots,\sE_n$ are vector bundles of rank $d$
on $S$ and we have morphisms
$$f_i:\sE_i \to \sE_{i+1},\quad f^i:\sE_{i+1} \to \sE_i$$
for each $i=1,\dots,n-1$. Given $s \in \Gamma(S,\sO_S)$, we say that
$\sE_{\bullet}=(\sE_i,f_i,f^i)_i$ is an {\bf $s$-linked chain} if the
following conditions are satisfied:
\begin{Ilist}
\itm For each $i=1,\dots,n$,
$$f_i \circ f^i = s \cdot \id, \text{ and } f^i \circ f_i = s \cdot \id.$$
\itm On the fibers of the $\sE_i$ at any point with $s=0$, we have that for
each $i=1,\dots,n-1$,
$$\ker f^i = \im f_i, \text{ and } \ker f_i = \im f^i.$$
\itm On the fibers of the $\sE_i$ at any point with $s=0$, we have that for
each $i=1,\dots,n-2$,
$$\im f_i \cap \ker f_{i+1}=(0),\text{ and }\im f^{i+1} \cap \ker f^i = (0).$$
\end{Ilist}
\end{defn}

We then define linked determinantal loci as follows.

\begin{defn}\label{defn:link-det} Let $\sE_{\bullet}$ be an $s$-linked
chain on a scheme $S$. Given $r,r_1,r_n>0$, suppose $\sF_1,\sF_n$ are 
vector bundles of rank $r_1$ and $r_n$ respectively, and let 
$g_1:\sE_1 \to \sF_n$ and $g_n:\sE_n \to \sF_n$ be any morphisms.
Then the associated \textbf{linked determinantal locus} is the closed 
subscheme of $S$ on
which the induced morphisms
\begin{equation}\label{eq:link-det}
\sE_i \to \sF_1 \oplus \sF_n
\end{equation}
have rank less than or equal to $r$ for all $i=1,\dots,n$.
\end{defn}

Note that the case $n=1$ and $\sF_n=0$ recovers the usual notion of
determinantal locus. It appears \textit{a priori} that the codimension
could be much larger than for a single determinantal locus, but it
turns out that the different determinantal conditions we are imposing 
are highly dependent. The case considered in Appendix A of \cite{os25} 
amounted to requiring that $g_1$ and $g_n$ be quotient maps onto bundles
of rank $r$. In that situation, it was shown that the classical codimension
bound also applies to linked determinantal loci. However, Cohen-Macaulayness 
does not appear to be addressable by the methods used in \textit{loc. cit}. 
Our main theorem is thus the following:

\begin{thm}\label{thm:main} If $S$ is Noetherian, then in the situation
of Definition \ref{defn:link-det}, every irreducible component of the
linked determinantal locus has codimension at most $(d-r)(r_1+r_n-r)$.
Moreover, if equality holds, and $S$ is Cohen-Macaulay then the linked 
determinantal locus is also Cohen-Macaulay.
\end{thm}

The proof proceeds by considering a suitable universal version of the
linked determinantal locus, and showing in essence that it is a 
``partial initial degeneration'' of the classical universal determinant 
locus. That is to say, in terms of ideals of minors we show that our
locus is defined by zeroing out monomials in each minor according to
a certain pattern, and that the resulting initial ideal is the same as
in the classical case. In comparison with other previously studied
variants of determinantal ideals, this appears to yield a rather distinct
direction of generalization.

As mentioned previously, using the construction given in \cite{os25},
we conclude from Theorem \ref{thm:main} that spaces of limit linear series
are flat and Cohen-Macaulay; see Theorem \ref{thm:lls-main} below.
We then conclude in Corollary \ref{cor:compare} that the two scheme
structures on spaces of limit linear series agree under typical 
circumstances (specifically, when the space has the expected dimension, 
the open subset of refined limit linear series is dense, and the
Eisenbud-Harris scheme structure is reduced). Finally, for the arguments
of \cite{c-l-t1} and \cite{c-l-p-t1}, the crucial consequence is
Corollary \ref{cor:proper-family}. This says that under the same
conditions as Corollary \ref{cor:compare}, given a one-parameter smoothing
of a curve $X_0$ of compact type, with generic fiber $X_{\eta}$, we have 
a flat, proper moduli space whose special fiber is the Eisenbud-Harris 
moduli space of limit linear series on $X_0$, and whose generic fiber is 
the usual linear series moduli space on $X_{\eta}$.

\section{Linked determinantal loci}

We consider the following situation:

\begin{sit} Let $R$ be a ring, and $s \in R$. Given $d,n,r_1,r_n$, and
nondecreasing integers $c_1 \leq c_2 \leq \dots c_{n-1}$, let
$R'$ be the polynomial ring over $R$ in variables $x_{i,j}$ with
$1 \leq i \leq r_1+r_n$ and $1 \leq j \leq d$. 
\end{sit}

\begin{defn} 
For $\ell=1,\dots,n$,
let $A_{\ell}$ be the $(r_1+r_n) \times d$ matrix over $R'$ defined by
$$\left(A_{\ell}\right)_{i,j} = 
\begin{cases} s^{e_{1,j,\ell}} x_{i,j}: & i \leq r_1 \\
s^{e_{1,j,\ell}} x_{i,j}: & i > r_1,\end{cases}$$
where 
$$e_{1,j,\ell}= \#\{m<\ell: j > d- c_m\}, \quad \text{ and } \quad
e_{2,j,\ell}= \#\{m\geq \ell: j \leq d- c_m\}.$$

Given $r>0$, let $I_r$ be the ideal of $R'$ defined by all $(r+1)\times(r+1)$
minors of all the $A_{\ell}$. In addition, let $J_r$ be the universal
determinantal ideal obtained from the $(r+1)\times(r+1)$ minors of the
matrix with $i,j$ entry equal to $x_{i,j}$.
\end{defn}

\begin{ex}\label{ex:basic} Consider the case $d=5$, $n=4$, $c_1=1$, 
$c_2=3$, $c_3=3$, $r=2$, $r_1=1$, $r_n=2$, and $s=0$. Then the matrices 
$A_{\ell}$ are as follows:

\begin{align*} & A_1=\begin{bmatrix} x_{1,1} & x_{1,2} & x_{1,3} & x_{1,4} & x_{1,5}\\
0 & 0 & 0 & 0 & x_{2,5} \\
0 & 0 & 0 & 0 & x_{2,5} \end{bmatrix}, \quad &
A_2=\begin{bmatrix} x_{1,1} & x_{1,2} & x_{1,3} & x_{1,4}&  0 \\
0 & 0 & x_{2,3} & x_{2,4} & x_{2,5} \\
0 & 0 & x_{3,3} & x_{3,4} & x_{3,5} \end{bmatrix}, \\
& A_3=\begin{bmatrix} x_{1,1} & x_{1,2} & 0 & 0 & 0 \\
0 & 0 & x_{2,3} & x_{2,4} & x_{2,5} \\
0 & 0 & x_{3,3} & x_{3,4} & x_{3,5} \end{bmatrix}, 
\quad & A_4=\begin{bmatrix} x_{1,1} & x_{1,2} & 0 & 0 & 0 \\
x_{2,1} & x_{2,2} & x_{2,3} & x_{2,4} & x_{2,5} \\
x_{3,1} & x_{3,2} & x_{3,3} & x_{3,4} & x_{3,5} \end{bmatrix}.\end{align*}

We do not list every minor from each of the $A_{\ell}$, but rather examine 
a selection of them with representative behavior.
The minors from the first three columns of the $A_{\ell}$ are 
$$0, \quad 0, \quad 0,\:\;\text{and}\:\; 
x_{1,1} x_{2,2} x_{3,3}-x_{1,1} x_{3,2} x_{2,3}-x_{2,1} x_{1,2} x_{3,3}
+ x_{3,1} x_{1,2} x_{2,3},$$ 
respectively.
The minors from the second, third and fourth columns of the 
$A_{\ell}$ are 
\begin{multline*} 0, \quad x_{1,2} x_{2,3} x_{3,4}-x_{1,2} x_{3,3} x_{2,4}, 
\quad x_{1,2} x_{2,3} x_{3,4}-x_{1,2} x_{3,3} x_{2,4}, \\
\:\;\text{and}\:\;x_{1,2} x_{2,3} x_{3,4}-x_{1,2} x_{3,3} x_{2,4},
\end{multline*} 
respectively.
Finally, the minors from the last three columns of the $A_{\ell}$ are 
$$0, \quad 
x_{1,3} x_{2,4} x_{3,5}-x_{1,3} x_{3,4} x_{2,5}-x_{2,3} x_{1,4} x_{3,5}
+ x_{3,3} x_{1,4} x_{2,5}, \quad 0,\:\;\text{and}\:\;0,$$ 
respectively.

The key aspects of these minors are the following:
\begin{itemize}
\item For each fixed minor position, the corresponding minors of the
different $A_{\ell}$ are always equal to zero or a unique nonzero value.
\item The unique nonzero value always occurs for some $\ell$, and always
contains the ``main diagonal'' term of the usual universal minor.
\item There need not be a unique choice of $\ell$ achieving the unique 
nonzero minor value, and no single value of $\ell$ generates nonzero
minors in all positions.
\end{itemize}
Note that the uniqueness of the nonzero minor value occurs in spite of the 
fact that a given nonzero minor may be generated by different patterns of 
zeroing out entries (as occurs with $A_2, A_3$ and $A_4$ for the minor
from the second, third and fourth columns).
\end{ex}

Throughout this section, we will use the lexicographic monomial order,
induced by the ordering that $x_{i,j}$ comes before $x_{i',j'}$ if
$i < i'$ or $i = i'$ and $j<j'$. In order to prove Theorem
\ref{thm:main}, the main result is then the following:

\begin{thm}\label{thm:universal-field} 
In the case $R=k$ is a field, then 
the initial ideal of $I_r$ coincides with
the initial ideal of $J_r$, and in particular the Hilbert functions coincide.
If $s \neq 0$, then also $R'/I_r \cong R'/J_r$.
\end{thm}

We prove this by exhibiting $I_r$ as a flat degeneration of $J_r$, using
the following construction.

\begin{defn}\label{defn:degeneration-matrix} Suppose $R=k[t]$, and let
$B$ be the $(r_1+r_n) \times d$ matrix over $R'$ defined by 
$$(B)_{i,j}=\begin{cases} t^{\epsilon_{1,j}} x_{i,j}: & i \leq r_1 \\
t^{\epsilon_{2,j}} x_{i,j} : & i>r_1,\end{cases}$$
where 
$$\epsilon_{1,j}= \#\{m:j>d-c_m\}, \quad \text{ and } \quad
\epsilon_{2,j}= \#\{m: j \leq d-c_m\}.$$
\end{defn}

Thus, $\epsilon_{1,j}=e_{1,j,n}$ and $\epsilon_{2,j}=e_{2,j,1}$.

\begin{ex}\label{ex:second} In the context of Example \ref{ex:basic}, we
have 
$$B=\begin{bmatrix} 
x_{1,1} & x_{1,2} & t^2 x_{1,3} & t^2 x_{1,4} & t^3 x_{1,5} \\
t^3 x_{2,1} & t^3 x_{2,2} & t x_{2,3} & t x_{2,4} & x_{2,5} \\
t^3 x_{3,1} & t^3 x_{3,2} & t x_{3,3} & t x_{3,4} & x_{3,5} \end{bmatrix}.$$

If we consider the minor from the first three columns of $B$, we get
\begin{multline*} t^4 x_{1,1} x_{2,2} x_{3,3} - t^4 x_{1,1} x_{3,2} x_{2,3}
- t^4 x_{2,1} x_{1,2} x_{3,3} 
\\ + t^4 x_{3,1} x_{1,2} x_{2,3} 
+ t^8 x_{2,1} x_{3,2} x_{1,3} - t^8 x_{3,1} x_{2,2} x_{1,3}.
\end{multline*}
Observe that if we factor out $t^4$ from this and then set $t=0$ we
obtain the nonzero minor from $A_4$ of Example \ref{ex:basic}.

The minor from the second, third and fourth columns of $B$ is
\begin{multline*} t^2 x_{1,2} x_{2,3} x_{3,4} - t^2 x_{1,2} x_{3,3} x_{2,4} 
-t^6 x_{2,2} x_{1,3} x_{3,4} \\
+ t^6 x_{3,2} x_{1,3} x_{2,4}
+t^6 x_{2,2} x_{3,3} x_{1,4} - t^6 x_{3,2} x_{2,3} x_{1,4},
\end{multline*}
and again factoring out $t^2$ and setting $t=0$ recovers the
corresponding nonzero minors from $A_2$, $A_3$ and $A_4$.

Finally, the minor from the last three columns of $B$ is
\begin{multline*} t^3 x_{1,3} x_{2,4} x_{3,5} - t^3 x_{1,3} x_{3,4} x_{2,5} 
- t^3 x_{2,3} x_{1,4} x_{3,5} \\
+ t^3 x_{3,3} x_{1,4} x_{2,5}
+ t^5 x_{3,3} x_{2,4} x_{1,5} - t^5 x_{2,3} x_{3,4} x_{1,5},
\end{multline*}
and factoring out $t^3$ and setting $t=0$ yields the corresponding
nonzero minor from $A_2$.

The above, together with similar analysis of the remaining minors, says
that in this example the ideal $I_r$ is governed in a suitable sense by 
the single matrix $B$, and it follows that $I_r$ is a flat degeneration 
of the universal determinantal ideal $J_r$.
\end{ex}

The following key lemma shows that the behavior of Examples \ref{ex:basic} 
and \ref{ex:second} generalizes.

\begin{lem}\label{lem:universal-field}
In the case that $R=k$ is a field and $s=0$, then if we fix
$i_{\bullet}=(i_0,\dots,i_r)$ with
$1 \leq i_0 <i_1 < \dots < i_r \leq r_1+r_n$ and
$j_{\bullet}=(j_0,\dots,j_r)$ with
$1 \leq j_0 < j_1 < \dots < j_r \leq d$, let $g_{i_{\bullet},j_{\bullet}}$
be obtained from the $(i_{\bullet},j_{\bullet})$ minor of $B$ by factoring
out the largest possible power of $t$, and then setting $t=0$.

Then for all $\ell$, the $(i_{\bullet},j_{\bullet})$ minor of $A_{\ell}$
is either equal to $0$ or to $g_{i_{\bullet},j_{\bullet}}$, and the latter
occurs at least once. In addition, $g_{i_{\bullet},j_{\bullet}}$ contains
the monomial $x_{i_0,j_0} x_{i_1,j_1} \cdots x_{i_r,j_r}$.
\end{lem}

\begin{proof} For each $\ell$, let $A'_{\ell}$ be the $(r+1)\times(r+1)$ 
matrix obtained from $A_{\ell}$ by taking the $(i,j)$ entries with 
$i \in \{i_0,\dots,i_r\}$ and $j \in \{j_0,\dots,j_r\}$, so that the
minor we are considering is $\det A'_{\ell}$. Let $m_1$ be
the number of $w$ such that $i_w \leq r_1$, and $m_2$ the number of
$w$ such that $i_w>r_1$, so that $m_1+m_2=r+1$. Then for any $\ell$,
there are $p_{1,\ell}, p_{2,\ell} \geq 0$ so that the first $m_1$ rows of 
the matrix $A'_{\ell}$ end in $p_{1,\ell}$ zeroes, while the last $m_2$ rows 
begin with $p_{2,\ell}$ zeroes (reading left to right). Explicitly,
we have 
$$p_{1,\ell}=\#\{w:j_w>d-c_{\ell-1}\}, \quad 
p_{2,\ell}=\#\{w:j_w \leq d-c_{\ell}\},$$ 
where here we should use the convention that $c_0=0$ and $c_{n}=d$.
In particular, we have $p_{1,\ell+1}+p_{2,\ell}=r+1$ for $\ell=1,\dots,n-1$.
Now, clearly $\det A'_{\ell}$ can only be nonzero if $r+1-p_{1,\ell} \geq m_1$
and $r+1-p_{2,\ell} \geq m_2$, and in fact, since the individual terms 
appearing in the determinant are distinct monomials, there cannot be any 
cancellation, so the converse holds as well. Then $p_{1,1}=0$, so 
if we choose $\ell$ maximal with $r+1-p_{1,\ell} \geq m_1$, 
and if $\ell<n$, then we have $p_{2,\ell}=r+1-p_{1,\ell+1} < m_1=r+1-m_2$,
so in this case $\det A'_{\ell} \neq 0$.
On the other hand, if $r+1-p_{1,n} \geq m_1$, then $p_{2,n}=0 \leq r+1-m_2$, 
so again $\det A'_{n} \neq 0$.

Thus, it remains to show that
for any $\ell$ such that $r+1-p_{1,\ell} \geq m_1$ and 
$r+1-p_{2,\ell} \geq m_2$, we have
that $\det A'_{\ell}= g_{i_{\bullet},j_{\bullet}}$, and that
$g_{i_{\bullet},j_{\bullet}}$ contains
the monomial $x_{i_0,j_0} x_{i_1,j_1} \cdots x_{i_r,j_r}$. Now, both
$\det A'_{\ell}$ and $g_{i_{\bullet},j_{\bullet}}$ are obtained by
omitting some monomials from the determinant of the $(r+1)\times (r+1)$
matrix with entries given by $x_{i_w,j_v}$, so we can prove the desired
statements by explicitly identifying which monomials are included in each.
First, for $g_{i_{\bullet},j_{\bullet}}$ we observe that since in 
$B$ the powers of $t$ in the first $r_1$ rows are nondecreasing from left
to right, and in the last $r_n$ rows are nonincreasing from left to right,
it is clear that (within the rows and columns determined by 
$i_{\bullet},j_{\bullet}$) if we take the top $m_1$ entries from the first 
$m_1$ columns (and consequently the remaining $m_2$ entries from the last 
$m_2$ columns)
we will simultaneously minimize the power of $t$ coming from the top $m_1$
rows and the bottom $m_2$ rows. In particular, the ``diagonal'' term 
obtained from
$(i_0,j_0),(i_1,j_1),\dots,(i_r,j_r)$ achieves the minimal possible power
of $t$ in the relevant minor, so $g_{i_{\bullet},j_{\bullet}}$ contains
$x_{i_0,j_0} x_{i_1,j_1} \cdots x_{i_r,j_r}$, as claimed. Moreover, we 
see that in order for a general term coming from 
$(i_0,j_{\sigma(0)}),\dots,(i_r,j_{\sigma(r)})$ for 
$\sigma \in \Sym(\{0,\dots,r\})$ to have the minimal power of $t$, we must
have $\epsilon_{1,j_{\sigma(0)}},\dots,\epsilon_{1,j_{\sigma(m_1-1)}}$ 
equal (as an unordered set with repetitions) to 
$\epsilon_{1,j_0},\dots,\epsilon_{1,j_{m_1-1}}$, 
and similarly for 
$\epsilon_{2,j_{\sigma(m_1)}},\dots,\epsilon_{2,j_{\sigma(r)}}$ and 
$\epsilon_{2,j_{m_1}},\dots,\epsilon_{2,j_r}$.

Now, suppose that 
$\epsilon_{1,j_{\sigma(0)}},\dots,\epsilon_{1,j_{\sigma(m_1-1)}}$ is not
equal to $\epsilon_{1,j_0},\dots,\epsilon_{1,j_{m_1-1}}$. Then 
let $w \leq m_1-1$ be minimal such that the entry in column $j_w$ is taken
from one of the bottom $m_2$ rows, so that $\epsilon_{1,j_w}$ occurs
strictly fewer times in the sequence
$\epsilon_{1,j_{\sigma(0)}},\dots,\epsilon_{1,j_{\sigma(m_1-1)}}$ than in
$\epsilon_{1,j_0},\dots,\epsilon_{1,j_{m_1-1}}$. If 
$\epsilon_{1,j_w}=\epsilon_{1,j_{m_1-1}}$, we conclude from minimality of
$\epsilon_{1,j_0},\dots,\epsilon_{1,j_{m_1-1}}$ that for 
some $w' \leq m_1-1$, we have 
$$\epsilon_{1,j_{\sigma(w')}}>\epsilon_{1,j_{m_1-1}}.$$
On the other hand, if $\epsilon_{1,j_w}<\epsilon_{1,j_{m_1-1}}$, we have that
$$\epsilon_{2,j_w} > \epsilon_{2,j_{m_1-1}} \geq \epsilon_{2,j_{m_1}},$$
and because we have taken the entry in the $j_w$ column from the bottom
$m_2$ rows, we conclude that for some $w' \geq m_1$, we have
$$\epsilon_{2,j_{\sigma(w')}}=\epsilon_{2,j_w}>\epsilon_{2,j_{m_1}}.$$

To summarize, we obtain a non-minimal power of $t$ for a given term if and 
only if $\epsilon_{1,j_{\sigma(w)}} \leq \epsilon_{1,j_{m_1-1}}$ for all 
$w <m_1$, and
$\epsilon_{2,j_{\sigma(w)}} \leq \epsilon_{2,j_{m_1}}$ for all $w \geq m_1$.
Thus, if we set $a$ so that $d-c_a < j_{m_1} \leq d-c_{a-1}$, and $b$ so 
that $d-c_b<j_{m_1+1} \leq d-c_{b-1}$,
then we are saying simply that a given term appears in
$g_{i_{\bullet},j_{\bullet}}$ if and only if 
\begin{equation}\label{eq:ineqs-1}
j_{\sigma(w)} \leq d-c_{a-1} \text{ for } w < m_1, \text{ and } \quad
j_{\sigma(w)} > d-c_b \text{ for }w \geq m_1.
\end{equation}

Now we consider the monomials in $\det A'_{\ell}$
for $\ell$ such that $r+1-p_{1,\ell} \geq m_1$ and 
$r+1-p_{2,\ell} \geq m_2$. It is clear that a given monomial
from entries $(i_0,j_{\sigma(0)}),\dots,(i_r,j_{\sigma(r)})$ occurs in
$\det A'_{\ell}$ if and only if 
\begin{equation}\label{eq:ineqs-2}
j_{\sigma(w)} \leq d-c_{\ell-1} \text{ for } w < m_1, \text{ and } \quad
j_{\sigma(w)} > d-c_{\ell} \text{ for }w \geq m_1.
\end{equation}
To compare \eqref{eq:ineqs-1} to \eqref{eq:ineqs-2}, we note that
$r+1-p_{1,\ell} \geq m_1$ is equivalent to saying that 
$j_{m_1} \leq d-c_{\ell-1}$,
and $r+1-p_{2,\ell} \geq m_2$ is the same as $p_{2,\ell} \leq m_1$, which
is the same as $j_{m_1+1}>d-c_{\ell}$,
so we can conclude that $\ell -1 \leq a-1$ and $\ell \geq b$, which is
to say that the relevant range for $\ell$ is $b \leq \ell \leq a$.
It then follows immediately that \eqref{eq:ineqs-1} implies 
\eqref{eq:ineqs-2}, and we wish to verify the converse.

If $a=b$, the converse is likewise immediate. However, we see from the
definitions that if $a >b$, then $d-c_b=d-c_{a-1}>d-c_a$, and then the
first part of \eqref{eq:ineqs-1} is equivalent to the second part.
In addition, we must have either $\ell=a$ or $c_{\ell}=c_b$. In the first 
case, we have that $c_{a-1}=c_{\ell-1}$, so the first part of 
\eqref{eq:ineqs-2} implies the first part of \eqref{eq:ineqs-1}, which
then implies the second part of \eqref{eq:ineqs-1} as well. But in the
second case, the second part of \eqref{eq:ineqs-2} implies the second
part of \eqref{eq:ineqs-1}, which then implies the first part as well.
\end{proof}

\begin{proof}[Proof of Theorem \ref{thm:universal-field}]
First suppose that $s \neq 0$. We see that in this case,
for $\ell <n$ the matrix $A_{\ell+1}$ may be obtained from $A_{\ell}$ 
by multiplying the righthand $c_{\ell}$ columns by $s$,
and dividing the bottom $r_n$ rows by $s$.
Thus, any given minor of $A_{\ell}$ is a power of $s$ times the 
corresponding minor of $A_{\ell+1}$, and we concude that $I_r$ is simply
equal to the ideal generated by the $(r+1)\times (r+1)$ minors of $A_1$.
As this is obtained from $J_r$ by rescaling the variables by powers of $s$,
we find that $R'/I_r \cong R'/J_r$, and also that
the initial ideals and Hilbert functions coincide.

Now, suppose that $s=0$. According to Lemma \ref{lem:universal-field},
in this case $I_r$ is generated by the $g_{i_{\bullet},j_{\bullet}}$, and
the initial terms of the latter agree with the initial terms of the
minors generating $J_r$. Now, let $J'_r$ be the flat degeneration 
of $J_r$ defined as in \S 15.8 of \cite{ei1} by the weight function 
assigning $\epsilon_{1,j}$ to the variable $x_{i,j}$ if $i \leq r_1$ and
$\epsilon_{2,j}$ to the variable $x_{i,j}$ if $i > r_1$. Then, it is clear 
from the definitions of the $g_{i_{\bullet},j_{\bullet}}$ and of $J'_r$
that $g_{i_{\bullet},j_{\bullet}} \in J'_r$ for each 
$(i_{\bullet},j_{\bullet})$. But by Theorem 15.17 and Exercise 20.14 of
\cite{ei1}, we have that the Hilbert functions of $J_r$ and of $J'_r$ 
coincide, so we conclude that the Hilbert function of $I_r$ is less than
or equal to the Hilbert function of $J_r$. On the other hand, the universal
minors form a Grobner basis for $J_r$ (see Theorem 5.3 of \cite{b-c1}), 
and their initial terms agree with those of the
$g_{i_{\bullet},j_{\bullet}}$, so we conclude that 
$\init J_r \subseteq \init I_r$. Again using invariance of Hilbert functions
under flat degenerations, we conclude that the Hilbert function of $I_r$
is greater than or equal to that of $J_r$, so they must be equal, and 
then we also have $\init J_r = \init I_r$, as desired.
\end{proof}

Theorem \ref{thm:main} then follows by standard reductions to known results
on the initial ideal of universal determinantal ideals. Indeed, we first
conclude:

\begin{cor}\label{cor:universal-field}
In the case $R=k$ is a field, then $R'/I_r$ is reduced and Cohen-Macaulay,
with codimension $(d-r)(r_1+r_n-r)$ in $R'$. If further $s \neq 0$ then
$R'/I_r$ is integral.
\end{cor}

\begin{proof} It is well known that $R'/\init J_r$ is reduced
and Cohen-Macaulay, of dimension
$(d-r)(r_1+r_n-r)$; see for instance Theorems 1.10, 5.3 and 6.7 of 
\cite{b-c1}. 
We conclude the same statements for $R'/I_r$ by Theorem 
\ref{thm:universal-field}, together with Proposition 3.12 of \cite{b-c1}.
In addition, for $s \neq 0$ we have $R'/I_r \cong R'/J_r$, and the latter
is integral (again by Theorem 1.10 of \cite{b-c1}).
\end{proof}

We next find:

\begin{cor}\label{cor:universal}
In the case that $R=\ZZ[t]$ and $s=t$, then $R'/I_r$ is flat
over $R$, and integral and Cohen-Macaulay, with codimension 
$(d-r)(r_1+r_n-r)$ in $R'$.
\end{cor}

\begin{proof} From Theorem \ref{thm:universal-field} we have
that all fibers of $R'/I_r$ over $R$ have the same Hilbert function. Since 
$R$ is reduced, it follows from Exercise 20.14 of \cite{ei1}
that $R'/I_r$ is flat over $R$. Given
flatness, the statements on irreducibility and codimension follows from
the corresponding statements on the generic fibers, which is a consequence
of the $s \neq 0$ statement of Corollary \ref{cor:universal-field}.
Again using flatness, the statements on Cohen-Macaulyness and reducedness
follow from the corresponding statements on fibers (see the Corollaries to
Theorems 23.3 and 23.9 of \cite{ma1}), which is again Corollary 
\ref{cor:universal-field}.
\end{proof}

In order to deduce our main theorem from the universal case, we recall 
Lemma 2.3 of \cite{o-t1}.

\begin{lem}\label{lem:structure} Suppose that $\sE_{\bullet}$
is $s$-linked. Let $c_i=\rk f_i$ for $i=1,\dots,n-1$, and by convention
set $c_0=0$, $c_n=d$. Also, for $j \leq i$ set 
$f_{j,i}= f_{i-1} \circ f_{i-2} \circ \dots \circ f_j$, and for $j \geq i$
set $f^{j,i}=f^{i} \circ f^{i+1} \circ \dots \circ f^{j-1}$.
Then locally on $S$, for $i=1,\dots,n$ there exist subbundles
$\sW_i \subseteq \sE_i$ of rank $c_i-c_{i-1}$ such that:
\begin{ilist}
\itm For $i=2,\dots,n-1$ we have that
$$\sW_i \cap \spn(\ker f_i, \ker f^{i-1})=(0),$$
and similarly $\sW_1 \cap \ker f_1=(0), \sW_n \cap \ker f^{n-1}=(0)$.
\itm For all $j<i$, the restriction of $f_{j,i}$ to $\sW_j$ is an isomorphism
onto a subbundle of $\sE_i$, and for $j>i$ the restriction of $f^{j,i}$ to
$\sW_j$ is an isomorphism onto a subbundle of $\sE_i$.
\itm The natural map
$$\left(\bigoplus_{j=1}^{i} f_{j,i}(\sW_j)\right) 
\oplus \left(\bigoplus_{j=i+1}^{n} f^{j,i}(\sW_j)\right) \to \sE_i$$
is an isomorphism for each $i$.
\end{ilist}
\end{lem}

\begin{proof}[Proof of Theorem \ref{thm:main}] The statement is local
on $S$, so we may assume that we have 
$\sW_i \subseteq \sE_i$ as in Lemma \ref{lem:structure}, and further
that $S=\Spec R$ is affine and the $\sW_i$ are free. Choose bases of 
the $\sW_i$ and use them to induce bases of the $\sE_i$ via Lemma
\ref{lem:structure} (iii), but with reversed ordering (so that basis
elements from $\sW_n$ come first, and those from $\sW_1$ last).
Choosing arbitrary bases of $\sF_1$ and $\sF_n$, 
we then have that the induced maps
$$\sE_i \to \sF_1 \oplus \sF_n$$
are given by matrices of the form of our $A_{\ell}$. These can be viewed
as induced by a map $\ZZ[\{x_{i,j}\},t] \to R$, and our linked determinantal
locus is then the pullback of the universal one under the corresponding
morphism $S \to \Spec \ZZ[\{x_{i,j}\},t]$. The theorem then follows
from Corollary \ref{cor:universal} by Theorem 3 and 
Proposition 4 of \cite{e-n1} (see also the introduction of \cite{e-h9}).
\end{proof}

We conclude with a couple of examples showing that the definition of linked
determinantal locus is somewhat delicate, in that minor variations will
invalidate the conclusion of Theorem \ref{thm:universal-field}.

\begin{ex}\label{ex:bad-zeroing} We first observe that a very similar 
pattern of zeroing out entries in a sequence of matrices can violate the
uniqueness of nonzero minors proved in Lemma \ref{lem:universal-field}.
Indeed, in Example \ref{ex:basic}, if we set $(A_3)_{1,3}$ to be $x_{1,3}$
instead of $0$, the minor from the last three columns is 
$x_{1,3} x_{2,4} x_{3,5} - x_{1,3} x_{3,4} x_{2,5}$, which does not
agree with the corresponding minor from $A_2$. Moreover, taking the 
difference yields
$x_{2,3} x_{1,4} x_{3,5} - x_{3,3} x_{1,4} x_{2,5}$, and since neither
of these monomials is in the initial ideal of $J_r$, we see that the 
conclusion of Theorem \ref{thm:universal-field} is also violated in this 
case.
\end{ex}

\begin{ex}\label{ex:bad-zeroing-2} We also see that if we zero out monomials
from the generators of the standard universal determinantal ideal, even if
the initial terms of each generator remain unchanged, in general the initial
ideals can change. For instance, if we let $J$ be the ideal generated by
the $2 \times 2$ minors of the matrix
$$A=\begin{bmatrix} x_{1,1} & x_{1,2} & x_{1,3} \\
x_{2,1} & x_{2,2} & x_{2,3} \end{bmatrix},$$
and we let 
$$I=(x_{1,1} x_{2,2}, x_{1,2} x_{2,3}, x_{1,1} x_{2,3}-x_{2,1} x_{1,3}),$$
then the initial terms of the generators are the same, and since we know the
minors are a Grobner basis for $J$, we conclude that 
$$\init I \supseteq \init J
=(x_{1,1} x_{2,2}, x_{1,2} x_{2,3}, x_{1,1} x_{2,3}).$$
However, in this case the containment is strict: we see that
$$x_{2,3} (x_{1,1} x_{2,2}) - x_{2,2} (x_{1,1} x_{2,3} - x_{2,1} x_{1,3})
= x_{2,1} x_{2,2} x_{1,3} \not \in \init J.$$
\end{ex}

\section{Applications to limit linear series}

We now apply our results to draw conclusions on spaces of limit linear
series. Because our results can be applied directly to limit linear
series moduli space constructions carried out in \cite{os25}, which are
explicitly in terms of linked determinantal loci, we have elected to keep
the presentation brief and not recall the rather lengthy definitions
leading up to the aforementioned constructions. However, we will below
briefly recall the Eisenbud-Harris limit linear series definition, so that
all objects relevant to our final conclusion, 
Corollary \ref{cor:proper-family}, have been defined.

Our main theorem deals with moduli spaces of limit linear series on
a one-parameter family $\pi:X \to B$ of curves somewhat more general than 
those of compact type. In the general setting, there will be some extra 
data denoted by $\bn$ and $(\sO_v)_v$, and an associated family
$\widetilde{\pi}:\widetilde{X} \to \widetilde{B}$, but in the compact type
case, these are irrelevant, and in particular we may assume 
$\widetilde{\pi}=\pi$. In either case, we denote the closed point of
$\widetilde{B}$ by $b_0$, and the special fiber of $\pi$ by $X_0$.
We will simultaneously treat the case of a single curve, where 
$\widetilde{B}=b_0$.

\begin{thm}\label{thm:lls-main} Suppose that we are in the situation of 
Theorem 6.1 of \cite{os25} (in particular, with $B$ the spectrum of a DVR), 
or of Theorem 5.9 of \cite{os25} (in which case $B$ is a point).
Suppose also that the space 
$G^r_{\bar{w}_0}(X_0,\bn,(\sO_v)_v)$ of limit linear series
on $X_0$ has the expected dimension $\rho$ at a given point $z$. Then 
the limit linear series moduli space 
$\widetilde{G}'^r_{\bar{w}_0}(\widetilde{X}/\widetilde{B},X_0,\bn,(\sO_v)_v)$
of Definition 6.3 of \cite{os25}
is Cohen-Macaulay at $z$, and flat over $\widetilde{B}$.
\end{thm}

In fact, the theorem applies also to higher-dimensional base schemes --
see Remark \ref{rem:higher-base}.

\begin{proof} According to Proposition 6.4 of \cite{os25}, 
$G':=
\widetilde{G}'^r_{\bar{w}_0}(\widetilde{X}/\widetilde{B},X_0,\bn,(\sO_v)_v)$
is described by the construction in the proof of Theorem 6.1 of \cite{os25}. 
This construction proceeds by constructing $G'$
as a closed subscheme of a scheme $G$ which is smooth over $\widetilde{B}$.
Furthermore, $\widetilde{B}$ is regular, so $G$ is likewise regular. 
The construction is given as an intersection of $(r+1)\deg D$ local
equations ensuring vanishing along an auxiliary divisor $D$, together with
$|V(\Gamma)-1|$ linked determinantal loci, each of expected codimension
$(r+1)(d+\deg D+1-g-(r+1))$. The numbers work out that in order for
$G^r_{\bar{w}_0}(X_0,\bn,(\sO_v)_v)$ to have dimension $\rho$ at $z$,
the aforementioned conditions must intersect with maximal codimension at $z$.
It thus follows that each individual condition cuts out a closed subscheme 
of $G$ of maximal codimension at $z$, and therefore by Theorem 
\ref{thm:main} we conclude that each of these closed subschemes is
Cohen-Macaulay at $z$. By Lemma 4.4 of \cite{o-h1} we conclude that the
intersection $G'$ is likewise Cohen-Macaulay at $z$.

Flatness is nontrivial only in the case that $B$ is positive-dimensional.
Then Theorem 6.1 of \cite{os25} implies that $G'$
is (universally) open over $B$ at $z$, so every irreducible component 
containing $z$ dominates $B$. Moreover, Cohen-Macaulayness implies that 
there are no imbedded components meeting $z$, so since $B$ is the spectrum 
of a DVR, we conclude that $G'$ is flat over $B$ at $z$, as desired.
\end{proof}

We now recall the Eisenbud-Harris definition in the compact type case.

\begin{defn} Given a curve $X_0$ of compact type, with dual graph
$\Gamma$, for $v \in V(\Gamma)$ let $Z_v$ denote the corresponding
component of $X_0$, and for $e \in E(\Gamma)$, let $P_e$ denote the
corresponding node.

Given $r,d>0$, a \textbf{limit linear series} (or more specifically, a 
\textbf{limit} $\fg^r_d$) on $X_0$ 
consists of a tuple $(\sL^v,V^v)_{v \in V(\Gamma)}$ of $\fg^r_d$s on the 
components $Z_v$, satisfying the condition that for each $e \in E(\Gamma)$
connecting vertices $v,v'$, we have:
\begin{equation}\label{eq:eh-ineq}
a^{(e,v)}_j+a^{(e,v')}_{r-j} \geq d \quad \text{ for } j=0,\dots,r,
\end{equation}
where $a^{(e,v)}$ denotes the vanishing sequence of $(\sL^v,V^v)$ at $P_e$.

We say the limit linear series is \textbf{refined} if \eqref{eq:eh-ineq}
is an equality for all $e \in E(\Gamma)$.
\end{defn}

In the compact type case, the Eisenbud-Harris definition of limit linear 
series leads to a natural scheme structure on the moduli space $G^r_d(X_0)$
of limit linear seres, 
as a union of closed subschemes ranging over all possible refined 
ramification conditions at the nodes. This definition has the advantage of
being very amenable to calculations, for instance in verifying reducedness.
The alternative definition introduced in \S 4 of \cite{os20}
and used in the statement of Theorem \ref{thm:lls-main} also
gives a scheme structure, which in principle could be different. The 
difficulty in comparing them arises from the union in the Eisenbud-Harris
case, as the functor of points of a union cannot be easily described. 
However, we can now conclude that under typical circumstances, the two 
scheme structures agree.

\begin{cor}\label{cor:compare} If $X_0$ is a curve of compact type,
suppose that we have the following conditions:
\begin{Ilist}
\itm $G^r_d(X_0)$ has the expected dimension $\rho$;
\itm the refined limit linear series are dense in $G^r_d(X_0)$;
\itm the Eisenbud-Harris scheme structure on $G^r_d(X_0)$ is reduced.
\end{Ilist}

Then the Eisenbud-Harris scheme structure on $G^r_d(X_0)$ coincides
with the scheme structure introduced in \cite{os20}.
\end{cor}

\begin{proof} According to Theorem \ref{thm:lls-main}, condition (I)
implies that the scheme structure of \cite{os20} is Cohen-Macaulay. Thus, 
in view of condition (III), it is
enough to show that the two scheme structures agree on a dense open subset,
since Cohen-Macaulayness will then imply reducedness. But Proposition 4.2.6
of \cite{os20} asserts that the two scheme structures agree on the refined
locus, so the desired result follows from condition (II).
\end{proof}

For applications such as \cite{c-l-t1} and \cite{c-l-p-t1}, the key point 
is the following immediate consequence of
Theorem \ref{thm:lls-main} and Corollary \ref{cor:compare}.

\begin{cor}\label{cor:proper-family} In the situation of Theorem 6.1 
of \cite{os25}, suppose further that the special fiber $X_0$ is of
compact type and the conditions of Corollary \ref{cor:compare} are
satisfied. Then
$\widetilde{G}'^r_{\bar{w}_0}(\widetilde{X}/\widetilde{B},X_0,\bn,(\sO_v)_v)$
is flat and proper over $B$, with special fiber equal to the
Eisenbud-Harris scheme structure on $G^r_d(X_0)$, and generic fiber 
equal to the classical space $G^r_d(X_{\eta})$ of linear series on $X_{\eta}$.
\end{cor}

\begin{rem}\label{rem:higher-base}
The moduli space construction of \cite{os25} restricted to the
case that $B$ was the spectrum of a DVR partly to avoid developing 
technical hypotheses such as the ``almost local'' condition of \S 2.2 of
\cite{os20}, and partly because in the non-compact-type case, it is 
necessary to impose more stringent conditions on the family to ensure that
every component of the curve comes from a divisor in the total family.
With that said, under suitable hypotheses the construction does generalize 
to higher-dimensional base schemes $B$, and in this case, Theorem
\ref{thm:lls-main} also generalizes. The only place where we used that
$B$ was the spectrum of a DVR was in arguing flatness, but as long as $B$
is regular this argument can be replaced via the use of 
Theorem 14.2.1 of \cite{ega43} and Proposition 6.1.5 of \cite{ega42}.

In the compact type context, it also suffices to assume that our family of 
curves has a single section (as opposed to a section through each component
of the special fiber, as assumed in \cite{os25}). The desired result being 
etale local, we can then
use an etale base change to produce sections through all components of
the special fiber.
\end{rem}

\bibliographystyle{amsalpha}
\bibliography{gen}

\newcommand{\noopsort}[1]{} \newcommand{\printfirst}[2]{#1}
  \newcommand{\singleletter}[1]{#1} \newcommand{\switchargs}[2]{#2#1}
\providecommand{\bysame}{\leavevmode\hbox to3em{\hrulefill}\thinspace}
\providecommand{\MR}{\relax\ifhmode\unskip\space\fi MR }
\providecommand{\MRhref}[2]{%
  \href{http://www.ams.org/mathscinet-getitem?mr=#1}{#2}
}
\providecommand{\href}[2]{#2}
\begin{thebibliography}{Oss14b}

\bibitem[BC03]{b-c1}
Winfried Bruns and Aldo Conca, \emph{Grobner bases and determinantal ideals},
  Commutative algebra, singularities and computer algebra (Siniai 2002), 2003,
  pp.~9--66.

\bibitem[CLPT]{c-l-p-t1}
Melody Chan, Alberto {L\'opez Martin}, Nathan Pflueger, and Montserrat
  {Teixidor i Bigas}, \emph{Genus of the {B}rill-{N}oether curve}, in
  preparation.

\bibitem[CLT]{c-l-t1}
Abel Castorena, Alberto {L\'opez Martin}, and Montserrat {Teixidor i Bigas},
  \emph{Invariants of the {B}rill-{N}oether curve}, preprint.

\bibitem[EH86]{e-h1}
David Eisenbud and Joe Harris, \emph{Limit linear series: Basic theory},
  Inventiones Mathematicae \textbf{85} (1986), no.~2, 337--371.

\bibitem[EH87a]{e-h5}
\bysame, \emph{Existence, decomposition, and limits of certain {W}eierstrass
  points}, Inventiones Mathematicae \textbf{87} (1987), 495--515.

\bibitem[EH87b]{e-h6}
\bysame, \emph{The {K}odaira dimension of the moduli space of curves of genus
  $\geq 23$}, Inventiones Mathematicae \textbf{90} (1987), 359--387.

\bibitem[Eis95]{ei1}
David Eisenbud, \emph{Commutative algebra with a view toward algebraic
  geometry}, Graduate Texts in Mathematics, vol. 150, Springer-Verlag, 1995.

\bibitem[EN67]{e-n1}
John Eagon and Douglas Northcott, \emph{Generically acyclic complexes and
  generically perfect ideals}, Proceedings of the Royal Society. Series A
  \textbf{299} (1967), 147--172.

\bibitem[GD65]{ega42}
Alexander Grothendieck and Jean Dieudonn\'e, \emph{{\'E}l\'ements de
  g\'eom\'etrie alg\'ebrique: {IV.} \'{E}tude locale des sch\'emas et des
  morphismes de sch\'emas, seconde partie}, Publications math\'ematiques de
  l'I.H.\'E.S., vol.~24, Institut des Hautes \'Etudes Scientifiques, 1965.

\bibitem[GD66]{ega43}
\bysame, \emph{{\'E}l\'ements de g\'eom\'etrie alg\'ebrique: {IV.} \'{E}tude
  locale des sch\'emas et des morphismes de sch\'emas, troisi\`eme partie},
  Publications math\'ematiques de l'I.H.\'E.S., vol.~28, Institut des Hautes
  \'Etudes Scientifiques, 1966.

\bibitem[HE71]{e-h9}
Melvin Hochster and John Eagon, \emph{{C}ohen-{M}acaulay rings, invariant
  theory, and the generic perfection of determinantal loci}, American Journal
  of Mathematics \textbf{93} (1971), 1020--1058.

\bibitem[HO08]{o-h1}
David Helm and Brian Osserman, \emph{Flatness of the linked {G}rassmannian},
  Proceedings of the AMS \textbf{136} (2008), no.~10, 3383--3390.

\bibitem[Kho]{kh1}
Deepak Khosla, \emph{Moduli spaces of curves with linear series and the slope
  conjecture}, preprint.

\bibitem[Mat86]{ma1}
Hideyuki Matsumura, \emph{Commutative ring theory}, Cambridge University Press,
  1986.

\bibitem[Oss14a]{os25}
Brian Osserman, \emph{Limit linear series for curves not of compact type},
  preprint, 2014.

\bibitem[Oss14b]{os20}
\bysame, \emph{Limit linear series moduli stacks in higher rank}, preprint,
  2014.

\bibitem[OT14]{o-t1}
Brian Osserman and Montserrat {Teixidor i Bigas}, \emph{Linked alternating
  forms and linked symplectic {G}rassmannians}, International Mathematics
  Research Notices \textbf{2014} (2014), no.~3, 720--744.

\end{thebibliography}

\end{document}